\documentclass[12pt]{amsart}
\usepackage{amssymb,amsthm,amsmath}
\usepackage{mdwlist}
\usepackage[utf8]{inputenc}
\usepackage[noend]{algorithmic}

\theoremstyle{plain}
\newtheorem{thm}{Theorem}[section]
\newtheorem{example}[thm]{Example} 
\theoremstyle{definition}
\newtheorem{definition}[thm]{Definition}
\newtheorem{remark}{Remark}
\newtheorem{algo}{Algorithm}

\begin{document} 
\title{Modulation in Tetradic Harmony and its Role in Jazz}
 
\author{Octavio A. Agustín-Aquino}
\address{Instituto de Física y Matemáticas,
Universidad Tecnológica de la Mixteca,
Huajuapan de León, Oaxaca, México}
\email{octavioalberto@mixteco.utm.mx}
\author{Guerino Mazzola}
\address{School of Music,
University of Minnesota,
Minneapolis, MN, USA
}
\email{mazzola@umn.edu}
\thanks{This work was partially supported by a grant from the
\emph{Niels Hendrik Abel Board}.}
\date{12 September 2018}


\begin{abstract}
After a quick exposition of Mazzola's quantum modulation model for the
so-called triadic interpretation of the major scale within the equal temperament, we study the model for the tetradic interpretation of the same scale. It is known that tetrads are fundamental for jazz music, and some classical objects for this kind of music are recovered.
\end{abstract}

\maketitle


\section{Introduction}

Mazzola's quantum model of modulation has been successful to explain, for instance,
the harmonic structure of the first movement of Beethoven's \emph{Hammerklavier} sonata
\cite[Part VI, Chapter 28]{gM02}, and it also explains the importance of melodic, harmonic, and major diatonic scales as extremes of modulatory freedom \cite[Part VI, Section 27.1.5]{gM02}.
Recently it has been used to guide the extension of traditional harmony to the
20-tone microtonal equal temperament (\cite{GTLLPM17,dG99}).

The notion of tonality in Western music is inextricably linked to the concept
of \emph{cadence}, which is a sequence of chords that has the intention of
asserting tonality. Considering that chords are identified by their position within
the scale with Roman numerals from the first (I) to the seventh (VII),
prominent examples of cadences are the \emph{authentic} V--I and the
\emph{plagal} IV--I cadences \cite[Ch. 10]{KP00} and Riemann's ``grosse Cadenz''
I--IV--I--V--I \cite[p. 52]{hR74}.

Together with the concept of tonality, there is the one of \emph{modulation}, that
consists in the transit from a source tonality to a target tonality in a
parsimonious way.

We now summarize Mazzola's theory of modulation (\cite[Part VI, Chapter 27]{gM02}) in order to explain how to extend it from the classical case of tonality with triads to tetrads or seventh chords, which are very important for jazz. First, we define an \emph{scale} $E$ as a subset of $\mathbb{Z}_{12}$ where $k=|E|$ is its cardinality. A \emph{triadic interpretation of the scale} 
$E$ (or \emph{tonality} of $E$), denoted by $E^{(3)}$, is a set of specific $3$-subsets of $E$ such that they cover $E$.

Each member of the triadic interpretation is called a \emph{degree}, and takes its
name according its position within the transversal in the natural order of
$E$ inherited from $\mathbb{Z}_{12}$, and it is usually denoted with a Roman
numeral. This general definition is usually further restricted to take into
account standard musical practice. More specifically, if $E=\{x_{i}\}_{i=0}^{k-1}$
is the indexing of $E$ following its natural order, then the triadic interpretation
\[
E^{(3)} = \{\{x_{i},x_{i+2\bmod k},x_{i+4\bmod k}\}\}_{i=0}^{k-1}
\]
corresponds to stacked thirds on top of each member of the scale. This is the interpretation we will use from now on.

\begin{example}

The subset $C=\{0,2,4,5,7,9,11\}\subseteq \mathbb{Z}_{12}$ is the C major scale,
and thus its triadic interpretation is
\[
 \{\{0,4,7\},\{2,5,9\},\{4,7,11\},\{7,11,2\},\{9,0,4\},\{11,2,5\}\}.
\]

For example, $I_{C}=\{0,4,7\}$, $II_{C}=\{2,5,9\}$, and so on.
\end{example}

A \emph{transposition} in $\mathbb{Z}_{12}$ is the map
\[
T^{a}:\mathbb{Z}_{12}\to\mathbb{Z}_{12}:x\mapsto x+a,
\]
and, in general, a \emph{special affine symmetry} is a map
\[
T^{a}.v:\mathbb{Z}_{12}\to\mathbb{Z}_{12}:x\mapsto vx+a
\]
where $v\in\{1,-1\}$.

From this moment on, we will consider only the orbit of the C scale with
respect to the group of transpositions extended pointwise to sets, and that we will call $Dia^{(3)}$. We have
\[
Dia^{(3)} = \{T^{a}(C^{(3)}):a\in\mathbb{Z}_{12}\}.
\]
\begin{remark}
Observe that actually
\[
 Dia^{(3)} = \{g(C^{(3)}):g\in T^{\mathbb{Z}_{12}}.\{\pm 1\}\}
\]
because $C^{(3)}$ is invariant with respect to the inversion symmetry $T^{4}.-1$.
\end{remark}

\section{Cadences}

\begin{definition}
Given $Dia^{(3)}$, a set $P$ of cadence parameters, and a map $\kappa:Dia^{(3)}\to P$, we say that $\kappa$ is \emph{cadential} in a tonality $E^{(3)}\in Dia^{(3)}$ if the
fiber $\kappa^{-1}(\kappa(E^{(3)}))$ is the singleton $\{E^{(3)}\}$, i. e., $\kappa$
is injective in $E^{(3)}$. The map $\kappa$ is a \emph{cadence} if it is cadential in
every tonality $E^{(3)}\in Dia^{(3)}$.
\end{definition}

We will take $P$ as covers of subsets of $\mathbb{Z}_{12}$
by degrees\footnote{This is an special case of the concept of \emph{global composition} \cite[Section 13.2]{gM02}.} of $E^{(3)}$. Thus a cadence is
completely determined by the indices of the degrees of
its image, and these are called \emph{cadential sets}. A cadential set is
\emph{minimal} if it does not have proper cadential subset.

For the case of $Dia^{(3)}$, in \cite{gM85} the following minimal cadences are calculated
\begin{multline*}
J_{1} = \{II,III\},
J_{2}=\{II,V\},J_{3}=\{III,IV\},\\
J_{4}=\{IV,V\},
J_{5}=\{VII\}.
\end{multline*}

Thus cadential sets are minimal amounts of information from which we can infer the tonality of a given musical fragment.

\section{Modulation}

Arnold Sch\"{o}nberg, in his \emph{Harmonielehre}, divides the process of
modulation in three stages \cite[p. 186]{aS74}:

\begin{quote}
\begin{enumerate*}
\item {[...]} introduction of such (neutral) chords as will
permit the turn to the new tonality [...].
\item the actual modulatory part, that is, the modulatory chords [...];
\item the consolidation, with the cadence establishing the target tonality.
\end{enumerate*}
\end{quote}

The chords mentioned in the second stage are the \emph{pivots} (Wendepunkte), and a
great deal of the modulation problem comes from finding them. Schönberg's definition
is the point of departure for Mazzola's quantum modulation model, describing the
connection between the source and the target tonalities through an affine symmetry.
From this symmetry and the cadence that consolidates the new tonality, it defines a
modulation \emph{quantum}, analogous to the ones that carry the fundamental
forces of nature and which are also are governed by symmetries.

\begin{definition}
A \emph{modulation} between two tonalities $S^{(3)}$ and $T^{(3)}$ is a pair $(g,\mu)$, where $g:S\to T$ is an special affine symmetry (called \emph{modulator}) that
induces an isomorphism between $S^{(3)}$ and $T^{(3)}$, and where $\mu$ is
a minimal cadential set of $T^{(3)}$.
\end{definition}

It can be proved that $g\in T^{\mathbb{Z}_{12}}.\{-1,1\}$, i. e., it is an
element of the special affine group of $\mathbb{Z}_{12}$.

\begin{definition}
A modulation \emph{quantum} $Q\subseteq \mathbb{Z}_{12}$ of a modulation $(g,\mu)$
satisfies the following:
\begin{enumerate*}
\item the modulator $g$ is a symmetry of $Q$, i. e., $g(Q) = Q$.
\item the set $T\cap Q$ is rigid (i. e., its symmetry group is trivial). The 
triads of the target tonality contained in this set are called the pivots.
\item The quantum $Q$ is the smallest set with the aforementioned properties.
\end{enumerate*}
\end{definition}

It is very remarkable that the model recuperates the pivots proposed by Sch\"{o}nberg in his treatise \cite[p. 574]{gM02}.

\section{The triadic nerve of the major scale}

\begin{figure}
\vspace{1em}
\begin{center}
\setlength{\unitlength}{0.15cm}
\begin{picture}(60,10)

\put(-5,0){\line(1,0){70}}
\put(-5,10){\line(1,0){70}}
\put(0,0){\line(-1,1){5}}
\put(0,0){\line(1,1){10}}
\put(10,10){\line(1,-1){10}}
\put(20,0){\line(1,1){10}}
\put(30,10){\line(1,-1){10}}
\put(40,0){\line(1,1){10}}
\put(50,10){\line(1,-1){10}}
\put(60,0){\line(1,1){5}}

\put(-10,-2){$A$}
\put(-10,9){$B$}

\put(67,9){$A$}
\put(67,-2){$B$}

\multiput(-5,0)(0,2){4}{\line(0,1){1}}
\put(-5,8){\vector(0,1){2}}
\multiput(65,3)(0,2){4}{\line(0,1){1}}
\put(65,2){\vector(0,-1){2}}

\put(0,-5){I}
\put(19,-5){V}
\put(39,-5){II}
\put(58,-5){VI}
\put(8,12){III}
\put(27,12){VII}
\put(48,12){IV}
\end{picture}
\end{center}
\vspace{1em}
\caption{Construction of the \emph{Harmonisches Band} by identification of two
sides $AB$ and $BA$ of a rectangle.}
\label{F:Moebius}
\end{figure}
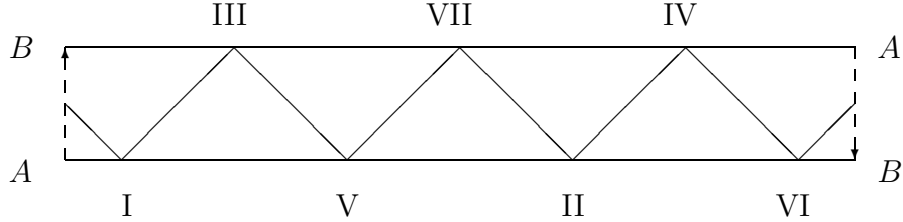

The \emph{nerve} \cite[pp. 128-129]{pA56} of $E^{(3)}$ is the simplicial complex of $E^{(3)}$ where every degree is a point and a set of points define
a simplex if their intersection is non-empty. Its geometrical realization (a topological space) turns out to be a Möbius strip
(Figure \ref{F:Moebius}), a fact discovered by Mazzola (\cite{gM85}) and implied from Schönberg's concept of \emph{Harmonisches Band}
\cite[p. 44]{aS74}. Furthermore, it is the triangulation with
the least number of vertices of the aforementioned space \cite[p. 32]{fC78}.
This topological construction has consequences for the function theory proposed by
Hugo Riemann, where he assigns to each cord its \emph{parallel} degree (Parallelklang) and its \emph{counterparallel}
(Leittonwechselklang), globally \cite[pp. 87-88]{dK06}. The unorientability of the Möbius strip makes such a global assignment impossible \cite[p. 324]{gM02}.

\section{The tetradic case}

In what follows, we extend Mazzola's model to seventh chords, that are tetrads of tones
separated by thirds. More specifically,
\[
E^{(4)} = \{\{x_{i},x_{i+2\bmod 7},x_{i+4\bmod 7},x_{i+6\bmod 7}\}\}_{i=0}^{6},
\]
which is crucial for jazz music \cite[p. 11,22]{RB05}. Cadential sets
are, in this case,
\begin{multline*}
J_{1} = \{I^{7},II^{7}\}, J_{2}=\{I^{7},IV^{7}\}, J_{3} = \{II^{7},III^{7}\},\\
J_{4} = \{III^{7},IV^{7}\}, J_{5} = \{V^{7}\}, J_{6} = \{VII^{7}\}.
\end{multline*}

\begin{remark}
The set of chords $\{II^{7},V^{7}\}$ is cadential, although not minimal. It is frequently used in jazz, like in \emph{Giant Steps} by John Coltrane \cite[p. 42]{RB05}.
\end{remark}

With the following algorithm (which is described in more
detail in \cite[p. 573]{gM02}) we can calculate all quanta.

\begin{algo}\hspace{1em}
\begin{algorithmic}[1]
\REQUIRE A modulation $(g,\mu)$ from the source tonality $E^{(3)}$ to the target tonality $F^{(3)}$.
\ENSURE The quantum $Q$ of the modulation $(g,\mu)$ if the modulation is quantized; a message of failure otherwise.
 \STATE $G:= \langle g \rangle, Q:= \emptyset$.
 \FORALL{$h\in G$}
   \STATE $Q:= Q\cup g\left(\bigcup\mu\right)$.
 \ENDFOR
 \STATE $M:= Q\cap F$.
 \STATE $\text{quantized}:= \texttt{true}$.
 \FORALL{$h\in T^{\mathbb{Z}_{12}}.\{-1,1\}\setminus\{1\}$}
  \IF{$h(M)=M$}
   \STATE $\text{quantized}:= \texttt{false}$.
  \ENDIF
 \ENDFOR
 \IF{$\text{quantized}$}
  \RETURN $Q$.
 \ELSE
  \PRINT The modulation is not quantized.
 \ENDIF
\end{algorithmic}
\end{algo}

\begin{table}
\begin{center}
\begin{tabular}{|c|c|c|c|c|}
\hline
Tr. & Cadence & Quantum & Modulator & Pivots\\
\hline
$1$ & $\{V^{7}\}$ & $\{0,2,3,5,6,8,9,11\}$ & $T^{5}.11$ & $\{III^{7},V^{7}\}$\\ 
\hline
$1$ & $\{VII^{7}\}$ & $\{0,2,3,5,6,7,10,11\}$ & $T^{5}.11$ & $\{VII^{7}\}$*\\ 
\hline
$2$ & $\{V^{7}\}$ & $\{1,2,4,5,7,9,11\}$ & $T^{6}.11$ & $\{II^{7},V^{7},VII^{7}\}$\\ 
\hline
$2$ & $\{VII^{7}\}$ & $\{1,2,4,5,7,11\}$ & $T^{6}.11$ & $\{II^{7},VII^{7}\}$\\ 
\hline
$3$ & $\{V^{7}\}$ &  $\{1,2,4,5,7,8,10,11\}$ & $T^{3}$ & $\{III^{7},V^{7}\}$\\ 
\hline
$3$ & $\{V^{7}\}$ &  $\{2,5,8,9,10,11\}$ & $T^{7}.11$ & $\{V^{7}\}$\\ 
\hline
$3$ & $\{VII^{7}\}$ &  $\{0,2,5,7,8,11\}$ & $T^{7}.11$ & $\{VII^{7}\}*$\\ 
\hline
$4$ & $\{V^{7}\}$ & $\{2,3,5,6,9,11\}$& $T^{8}.11$ & $\{V^{7}\}$\\ 
\hline
$5$ & $\{VII^{7}\}$ & $\{2,4,5,7,10,11\}$& $T^{9}.11$ & $\{II^{7},VII^{7}\}$\\  
\hline
$6$ & $\{V^{7}\}$ & $\{1,2,5,7,8,11\}$ &$T^{6}$& $\{V^{7}\}$\\ 
\hline
$6$ & $\{V^{7}\}$ & $\{1,2,5,8,9,11\}$ &$T^{10}.11$& $\{V^{7}\}$\\ 
\hline
$6$ & $\{VII^{7}\}$ & $\{2,3,5,8,9,11\}$ &$T^{6}$ & $\{VII^{7}\}$\\ 
\hline
$6$ & $\{VII^{7}\}$ & $\{2,3,5,7,8,11\}$ &$T^{10}.11$ & $\{VII^{7}\}$\\ 
\hline
$7$ & $\{V^{7}\}$ & $\{0,2,5,6,9,11\}$ &$T^{11}.11$ & $\{III^{7},V^{7}\}$\\ 
\hline
$8$ & $\{VII^{7}\}$ & $\{1,2,5,7,10,11\}$ & $T^{0}.11$ & $\{VII^{7}\}$\\ 
\hline
$9$ & $\{V^{7}\}$ & $\{1,2,4,5,7,8,10,11\}$ & $T^{9}$ & $\{III^{7},V^{7}\}$\\ 
\hline
$9$	& $\{V^{7}\}$ & $\{2,4,5,8,9,11\}$ & $T^{1}.11$	&$\{V^{7}\}*$\\ 
\hline
$9$ & $\{VII^{7}\}$	& $\{2,5,6,7,8,11\}$ &$T^{1}.11$ &$\{VII^{7}\}$\\ 
\hline
$10$ &	$\{V^{7}\}$ & $\{0,2,3,5,9,11\}$ & $	T^{2}.11$ & $	\{III^{7},V^{7}\}$\\ 
\hline
$10$ &	$\{VII^{7}\}$ & $\{0,2,3,5,7,9,11\}$ & $	T^{2}.11$ & $\{	III^{7},V^{7},VII^{7}\}$\\ 
\hline
$11$ &	$\{V^{7}\}$ & $\{1,2,4,5,6,9,10,11\}$ & $	T^{3}.11$ & $\{	V^{7}\}$*\\ 
\hline
$11$ &	$\{VII^{7}\}$ & $\{1,2,4,5,7,8,10,11\}$ & $	T^{3}.11$ & $\{II^{7},VII^{7}\}$\\ 
\hline
\end{tabular}
\end{center}
\caption{Summary of quantized modulations for the tetradic tonalities over the
major scale.}
\label{T:Modulaciones}
\end{table}

The output for each quantized modulation is shown in Table \ref{T:Modulaciones} and leads us to the following insights.

\begin{enumerate*}
\item All quantized modulations for the tetradic interpretations stem from the $\{V\}$
and $\{VII\}$ cadences.
\item In the theory as formulated in \cite{gM02}, it is required that $T\cap Q$
to be covered by triads, but this is not required in the original definition
of \cite{gM85}. In fact, the pivots marked with an asterisk in Table \ref{T:Modulaciones}
are not enough to cover $T\cap Q$.
\item The modulation quanta of modulations
\[
 (T^{\{3,5\}}.11,\{V^{7}\})\quad\text{and}\quad(T^{\{3,9\}},\{V^{7}\})
\]
are what in jazz are known as \emph{diminished scales}, which
are $8$-tone scales in which the notes ascend in alternating intervals of a whole step and a half step \cite[p. 31]{RB05}.
\item Modulations with cadences $\{V^{7}\}$ and $\{VII^{7}\}$ to a distance of tritone correspond to the jazz artifice of \emph{tritone substitution} \cite[p. 147]{dL93}. More specifically, the pivot of the target tonality is the tritone substitution of the same degree of the source tonality.
\item A particular case of a quick modulation towards a tonality a major second upwards known as \emph{chaining} \cite[p. 46]{RB05} is recovered as a quantized modulation for cadence $\{V^{7}\}$ with pivots $\{II^{7},V^{7}\}$.
\end{enumerate*}

\section{The tetradic nerve of the major scale}
If we calculate the nerve of $E$ with the tetradic cover $E^{(4)}$, we see that
all degrees intersect in at least one point, thus the $1$-skeleton of such a nerve is the complete graph $K_{7}$. This is a toroidal graph \cite[p. 71]{aG85}, and its embedding in such a surface (Figure \ref{F:Toro}) is the simplicial triangulation with the smallest number of vertices \cite[p. 277]{KMS15}. Moreover, there are seven $3$-dimensional simplices or tetrahedra
\begin{multline*}
S_{1}=\{I^{7},II^{7},IV^{7},VI^{7}\},
S_{2}=\{II^{7},III^{7},V^{7},VII^{7}\},\\
S_{3}=\{III^{7},IV^{7},VI^{7},I^{7}\},
S_{4}=\{IV^{7},V^{7},VII^{7},II^{7}\},\\
S_{5}=\{V^{7},VI^{7},I^{7},III^{7}\},
S_{6}=\{VI^{7},VII^{7},II^{7},IV^{7}\},\\
S_{7}=\{VII^{7},I^{7},III^{7},V^{7}\},
\end{multline*}
which can be arranged in a cycle $S_{5},S_{7},S_{2},S_{4},S_{6},S_{1},S_{3},S_{5}$ that forms a solid torus. Now, even although $K_{7}$ is embeddable in an orientable surface, it is not possible to formulate a Riemann global theory of chord functions: the obstacle is the fact that the $1$-skeleton graph of the nerve contains the cycle $I^{7},V^{7},IV^{7},VI^{7},III^{7},VII^{7},II^{7},I^{7}$, which makes impossible to define paralellism in a consistent way.

\section*{Acknowledgment}

We thank Jesús David Gómez Téllez for detecting some mistakes in a preliminary version of this article and suggesting some improvements.

\begin{figure}
\begin{center}
\setlength{\unitlength}{0.2cm}
\begin{picture}(30,40)

\put(0,0){\vector(0,1){40}}
\put(0,0){\vector(1,0){30}}
\put(0,40){\vector(1,0){30}}
\put(30,0){\vector(0,1){40}}

\put(0,10){\line(1,-1){10}}
\put(0,10){\line(2,-1){20}}
\put(0,10){\line(1,1){10}}
\put(0,30){\line(1,-1){10}}
\put(0,40){\line(1,-2){10}}

\put(10,40){\line(0,-1){20}}
\put(10,40){\line(1,-2){10}}

\put(10,20){\line(1,0){10}}
\put(10,20){\line(1,-2){10}}

\put(30,30){\line(-1,1){10}}
\put(30,30){\line(-2,1){20}}
\put(30,30){\line(-1,-1){10}}

\put(30,10){\line(-1,1){10}}
\put(30,0){\line(-1,2){10}}

\put(20,0){\line(0,1){20}}

\put(-3,-1){III}
\put(-3,9){IV}
\put(-3,29){II}
\put(-3,39){III}

\put(9,-2){VII}
\put(6,19.5){I}
\put(9,41){VII}

\put(19,-2){VI}
\put(22,19.5){V}
\put(19,41){VI}

\put(31,-1){III}
\put(31,9){IV}
\put(31,29){II}
\put(31,39){III}

\end{picture}
\end{center}
\caption{Construction of the embedding of $K_{7}$ in the torus by identifications
of the opposite sides of a rectangle.}
\label{F:Toro}
\end{figure}

\bibliographystyle{abbrv}
\bibliography{seminario_utm_2017}{}

\end{document}